\documentclass{jcmlatex}
\usepackage{graphicx}
\usepackage{algorithm,algorithmic} 
\usepackage{cleveref} 
\usepackage{colortbl}
\usepackage{subfigure}
\usepackage{listings} 
\usepackage{arydshln} 

\Crefname{assumption}{Assumption}{Assumptions}

\def\JCMvol{xx}
\def\JCMno{x}
\def\JCMyear{200x}
\def\JCMreceived{xxx}
\def\JCMrevised{xxx}
\def\JCMaccepted{xxx}

\begin{document}

\markboth{W. S. TENG AND Q. N. LI}{Template for Journal Of Computational
Mathematics}

\title{A MATRIX OPTIMIZATION METHOD FOR BLIND EXTRACTION OF EXTERNAL EQUITABLE PARTITIONS FROM LOW PASS GRAPH SIGNALS}

\author{Wenshun Teng
\thanks{School of Mathematics and Statistics, Beijing Institute of Technology, Beijing, China \\ Email: teng\_wenshun@163.com}
\and
Qingna Li\footnote{Corresponding author. The author's research is supported by NSFC 12071032 and NSFC 12271526}
\thanks{School of Mathematics and Statistics/Beijing Key Laboratory on MCAACI, Beijing Institute of Technology, Beijing, China\\ Email: qnl@bit.edu.cn}}

\maketitle

\begin{abstract}
Seeking the external equitable partitions (EEPs) of networks under unknown structures is an emerging problem in network analysis. The special structure of EEPs has found widespread applications in many fields such as cluster synchronization and consensus dynamics. While most literature focuses on utilizing the special structural properties of EEPs for network studies, there has been little work on the extraction of EEPs or their connection with graph signals. In this paper, we address the interesting connection between low pass graph signals and EEPs, which, as far as we know, is the first time.
We provide a method BE-EEPs for extracting EEPs from low pass graph signals and propose an optimization model, which is essentially a problem involving nonnegative orthogonality matrix decomposition. We derive theoretical error bounds for the performance of our proposed method under certain assumptions and apply three algorithms to solve the resulting model, including the K-means algorithm, the practical exact penalty method and the iterative Lagrangian approach. Numerical experiments verify the effectiveness of the proposed method. Under strong low pass graph signals, the iterative Lagrangian and K-means perform equally well, outperforming the exact penalty method. However, under complex weak low pass signals, all three perform equally well.
\end{abstract}



\section{Introduction}

The extraction of external equitable partitions (EEPs) is an emerging problem in network analysis \cite{num1}. EEPs are a graph topological structure, which is a generalization of the well-known equitable partitions (EPs) \cite{no27}, also known as almost equitable partitions (AEPs) \cite{2-3eep} in graph theory. Cardoso et al. \cite{2-3eep} introduced AEPs and provided their existence. Previous research has focused on various areas such as graph coloring \cite{num17} and automorphism \cite{num18} in graph theory. Recent studies have mostly focused on network dynamic processes. EEPs have an intriguing geometric interpretation, revealing spectral properties of the Laplacian matrix \cite{num4,num5}. This property opens up more possibilities for network research, such as cluster synchronization \cite{no11,no28}, consensus dynamics \cite{no12,no29,no30}, network control problems \cite{new-no11,num12,num13}, and epidemic spreading \cite{num14,num15}. EEPs have even been utilized to aid in studying the microstructure of molecules and crystalline solids \cite{num16}. Additionally, the research by Schaub et al. \cite{2-4eep} involves the task of finding approximate EEPs. To address this, they propose minimizing a projection error problem and utilize K-means clustering to solve it.

On the other hand, the low pass graph signal model is a hot topic in the field of graph signal processing (GSP) \cite{no17,no18}, with widespread applications in various domains such as economics, social networks, and power systems \cite{lowpass}. Many physical and social processes are naturally characterized by low pass graph filters, such as diffusion model \cite{diffusion}, opinion dynamics \cite{opinion}, and soon. Recent research mainly focuses on two aspects:
(1) Due to the significant role played by low pass graph filters in tasks related to sampling \cite{num19}, denoising \cite{num20}, graph topology inference \cite{num21}, and graph neural networks \cite{num22}, most studies typically develop relevant graph signal processing algorithms with low pass attributes as a default assumption.
(2) Considering scenarios where the properties of observed graph signals are unknown, such as systems under attack \cite{num23} or observed data being corrupted \cite{num24}, many researchers studied the blind detection problem of low pass graph signals. This involves detecting whether a set of observed graph signals is generated by a low pass graph filter, as studied in \cite{num25,num26,num27}. Our work is related to graph topology inference, with recent studies focusing on blind inference of communities \cite{no32,no31} and centrality \cite{no22,no24,IIR}.

Having summarized the existing work on EEPs and low pass graph signal models, it seems that little work has been done on the connection of EEPs and low pass graph signal models. A natural question is whether the two topics can be related to each other and what is the exact connection between the two? Scholkemper et al. \cite{no26} established a link between EPs and a specific graph signal model, which inspired our research. In this paper, we aim to answer this question by establishing the mathematical model of EEPs via low pass graph signals.

The contribution of this paper is as follows. Firstly, we reveal the interesting connection between EEPs and low pass graph signals. Secondly, we propose a method for learning the EEP via low pass graph signals, which involves constructing an optimization model. Finally, we provide theoretical analysis for the proposed method and apply three algorithms to solve it. This paper opens a door to investigate EEPs from a new point of view, and provide a new approach to extract the EEP of a graph.

The organization of the rest of this paper is as follows. In \Cref{sec2}, we introduce some preliminaries about EEPs of a graph and graph signal models. In \Cref{sec3}, we establish the optimization model for extraction of EEPs via low pass graph filters and provide a method, namely BE-EEPs, for extracting EEPs from low pass graph signals. In \Cref{sec4}, we derive theoretical error bounds for the performance of BE-EEPs under certain assumptions. In \Cref{sec5}, we discuss how to solve the proposed optimization model. In \Cref{sec6}, we present some numerical results. Final conclusions are made in \Cref{sec7}.

\section{Preliminaries}\label{sec2}
In this section, we introduce some preliminaries about EEPs of a graph and graph signal models.


\subsection{External Equitable Partitions of a Graph}\label{sec2-1}
\begin{figure}[htp]
\centering
\subfigure[The original graph $G$]
{
 	\begin{minipage}[b]{.5\linewidth}
        \centering
        \includegraphics[scale=0.8]{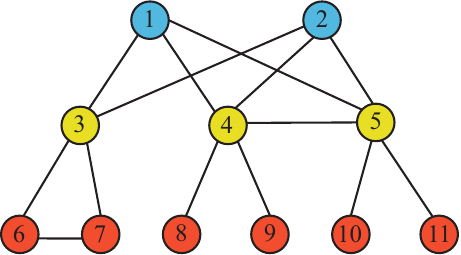}
    \end{minipage}
}
\subfigure[Quotient graph $G/\pi$]
{
 	\begin{minipage}[b]{.45\linewidth}
        \centering
        \includegraphics[scale=0.8]{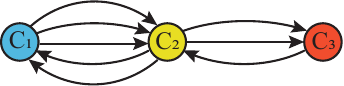}
    \end{minipage}
}
\caption{The original graph $G$ with $n=11$ vertices with an EEP $\pi$ into three cells (indicated with colors) and its associated quotient graph $G/\pi$.}
\label{fig-quotient}
\end{figure}

For an undirected simple graph $G=(V(G), E(G))$, $V(G)=\{1, 2, \cdots, n\}$ is a nonempty set of vertices, $E(G)=\{(i, j) \ |\ i, j=1, \cdots, n, i\ {\rm and }\ j \ {\rm is}\ {\rm connected} \}$ is a set of edges. A partition $\pi$ of $V(G)$ with cells $\mathbb{C}= \{C_{1}, C_{2}, \cdots, C_{r}\}$ means that $V(G)$ is divided into disjoint sets $C_{i}$, $i=1, \cdots, r$ such that $C_{i}\cap C_{j}=\emptyset$, $i\neq j$, $j=1, \cdots, r$ and $V(G)=\bigcup_{i} C_{i}$. We say that a partition $\pi$ of $V(G)$ with cells $C_{1}$, $\cdots$, $C_{r}$ is an $equitable$ $partition$ (EP) if the number of neighbours in $C_{j}$ of a vertex $u$ in $C_{i}$ is a constant $b_{ij}$, independent of $u$, where $i, j=1, \cdots, r$. Similarly, a partition $\pi$ of $V(G)$ with cells $C_{1}$, $\cdots$, $C_{r}$ is an $external$ $equitable$ $partition$ (EEP) if the number of neighbours in $C_{j}$ of a vertex $u$ in $C_{i}$ is a constant $b_{ij}$, independent of $u$, where $i\neq j$ and $i, j=1, \cdots, r$. Note that compared with EP, in EEP, only for $i\neq j$, it is required that all the vertices in the class $C_{i}$ have the same number of neighbors in $C_{j}$, $i, j=1, 2, \cdots, r$. Let $G/\pi$ denote the quotient graph in which each vertex represents a cell of the EEP $\pi$ and the arcs between the new vertices indicate the links between cells in the original graph $G$ (see \Cref{fig-quotient}). Note that the quotient graph $G/\pi$ is a multi-digraph.

To characterize a partition, a $binary$ $indicator$ $matrix$ \cite{bim} $H\in\mathbb{R}^{n\times r}$ is often used, where
\begin{equation}
\label{eqH}
H_{ij} =\left\{\begin{array}{ll}
		1,&{\rm if}\ i\in C_{j},\\
		0,&{\rm otherwise}.
		\end{array}	\right.
\end{equation}
Next, for an EEP $\pi$ of $G$, let $A^{G}\in \{0, 1\}^{n\times n}$ be the adjacency matrix of the graph $G$ in which $(A^{G})_{ij} = 1$ if $(i,j)\in E(G)$ and $(A^{G})_{ij} = 0$ otherwise. Let $A^{G/\pi}\in\mathbb{R}^{r\times r}$ be the adjacency matrix of the quotient graph $G/\pi$ in which $(A^{G/\pi})_{ij} = b_{ij}$ if $i\neq j$ and $(A^{G/\pi})_{ij} = 0$ otherwise, where $b_{ij}$ is the constant defined in EEP. Let $L^{G}\in\mathbb{R}^{n\times n}$ and $L^{G/\pi}\in\mathbb{R}^{r\times r}$ be the Laplacian matrix of the graph $G$ and the quotient graph $G/\pi$, that is,
$$L^{G}=D^{G}-A^{G}\ {\rm and} \ L^{G/\pi}=D^{G/\pi}-A^{G/\pi},$$
where $D^{G}={\rm Diag}({A^{G}{\rm\textbf{1}}})$ and $D^{G/\pi}={\rm Diag}({A^{G/\pi}{\rm\textbf{1}}})$ are diagonal matrices with the diagonal elements being $A^{G}{\rm\textbf{1}}$ and $A^{G/\pi}{\rm\textbf{1}}$ respectively, and ${\rm\textbf{1}}$ is an all-one vector of proper size. Although $L^{G}$ is symmetric, $L^{G/\pi}$ will be asymmetric in general. The Laplacian matrix of the quotient graph $G/\pi$ ignores the internal connectivity. Note that, from the definition of the Laplacian, there is always a trivial EEP in which the whole graph is grouped into one class, i.e., $H=\textbf{1}$ and $L^{G/\pi}=0$. Moreover, every EP is necessarily an EEP, while the converse is not true.

Based on \cite{no27, 2-3eep, 2-4eep}, we have the following properties which will be used later. More properties can be found in \cite{no27, 2-3eep, no11, 2-4eep}.
\begin{proposition}\label{prop-ep-eep}
Let $G$ be a graph, and $\pi$ be a partition of $V(G)$ encoded by the indicator matrix $H$.
The following results hold.
\begin{itemize}
\item[{\rm (i)}] The partition $\pi$ is an EEP if and only if $L^{G}H=HL^{G/\pi}$.
\item[{\rm (ii)}] If $\pi$ is an EEP, then there exist $r$ eigenvectors contained in $V^{G}=[v_{1},\cdots,v_{r}]\in\mathbb{R}^{n\times r}$ for the Laplacian matrix $L^{G}$ such that $V^{G}=HV^{G/\pi}$, where $V^{G/\pi}=[v_{1}^{G/\pi},\cdots,v_{r}^{G/\pi}]$ contains $r$ eigenvectors of $L^{G/\pi}$. 
\end{itemize}
\end{proposition}
\begin{definition}
For an EEP denoted as $\pi$, if there exists an eigenvector $v\in\mathbb{R}^{n}$ of $L^{G}$ and an eigenvector $v^{G/\pi}\in\mathbb{R}^{r}$ of $L^{G/\pi}$ such that $v=Hv^{G/\pi}$, then $v$ is referred to as a structural eigenvector of $L^{G}$.
\end{definition}

\subsection{Graph Signal Model}\label{sec2-2}
Recall the undirected graph $G=(V(G), E(G))$ defined in \Cref{sec2-1}. Since $L^{G}$ is real symmetric and positive semidefinite \cite{no32}, it can be written as:
\begin{equation}\label{2.2}
L^{G}=V\Lambda V^{\top},
\end{equation}
where $V=[v_{1}, v_{2}, \cdots, v_{n}]$ is the eigenvector matrix that contains the eigenvectors as columns, and $\Lambda={\rm Diag}[\lambda_{1}, \lambda_{2}, \cdots, \lambda_{n}]$ is the diagonal eigenvalue matrix. The eigenvalues are ordered as $0=\lambda_{1} \leq\cdots\leq\lambda_{n}$. Note that $V$ is an orthogonal matrix.

\begin{figure}[htp]
  \centering
  \includegraphics[width=0.5\textwidth]{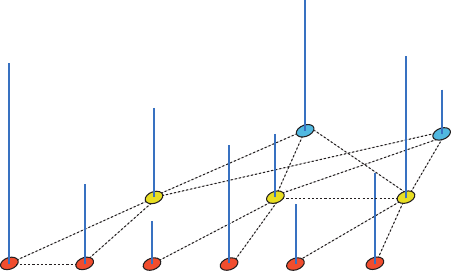}\\
  \caption{A graph signal on the vertices of graph $G$ in \Cref{fig-quotient}. The height of each blue bar represents the magnitude of the signal values.}
  \label{fig-G-signal}
\end{figure}

A signal on the graph $G$ is defined as a function on the vertices of $G$, $f: V(G)\rightarrow\mathbb{R}$. For the sake of description, a signal on the graph $G$ can be represented as a vector $\textbf{x}:=[x_{1}, x_{2}, \cdots, x_{n}]\in \mathbb{R}^{n}$ that attaches a value $x_{i}$ to the vertex $i$ in $V(G)$, $i=1, 2, \cdots, n$. The graph signal in \Cref{fig-G-signal} is one such example.
The output graph signal $\textbf{y}=[y_{1}, y_{2}, \cdots, y_{n}]\in \mathbb{R}^{n}$ of a graph filter $\mathcal{H}$ applied to the input graph signal $\textbf{x}$ is
\begin{equation}
\textbf{y}= \mathcal{H}\textbf{x},
\end{equation}
where $\mathcal{H}$ is represented by matrices called a $graph$ $shift$ $operator$\footnote{Graph shift operator is commonly used as the graph Laplacian matrix $L^{G}$ and the adjacency matrix $A^{G}$. See \cite{no17} for an overview on the subject.} \cite{no31}, which can be used to define linear graph filters.
Here, we use the graph Laplacian matrix $L^{G}$ as a graph shift operator.
A graph filter $\mathcal{H}$ of order $T$ ($T\leq+\infty$) can be expressed as a polynomial of matrix in the graph shift operator $L^{G}$ of degree $T$ (recall \eqref{2.2}):
\begin{equation}\label{H(L)}
\mathcal{H}(L^{G}):=\sum\limits^{T}_{t=0}a_{t}(L^{G})^{t}=VDiag(h(\lambda_{1}),\cdots,h(\lambda_{n}))V^{\top},
\end{equation}
where $h(\cdot):\ \mathbb{R}\rightarrow\mathbb{R}$ is the scalar generating polynomial defined by
\begin{equation}
h(\mu)=\sum\limits^{T}_{t=0}a_{t}\mu^{t},\ \mu\in\mathbb{R},
\end{equation}
and $a_{t}\in\mathbb{R}$, $t=0,\cdots,T$ is the filter weight. Let the filtered graph signals be given by $\textbf{y}=\mathcal{H}(L^{G})\textbf{x}$.
Now we describe the definition of a low pass graph filter.

\begin{definition}{\rm\cite{lowpass}}
\label{deflowpass}
We say that a graph filter $\mathcal{H}(L^{G})$ is $r$-low-pass if for any $1\leq r\leq n-1$,
\begin{equation}\label{eta}
\eta_{r}: =\frac{\max\{|h(\lambda_{r+1})|,\cdots,|h(\lambda_{n})|\}}{\min\{|h(\lambda_{1})|,\cdots,|h(\lambda_{r})|\}}\in\left.\left[0,1\right.\right).
\end{equation}
The ratio $\eta_{r}$ quantifies the $strength$ of the low pass graph filter. Generally, if $\eta_{r}\approx0$, then we can say the low pass graph filter is $strong$; if $\eta_{r}\approx1$, then the low pass graph filter is $weak$.
\end{definition}

\Cref{deflowpass} implies that the eigenvectors of a $r$-low-pass graph filter $\mathcal{H}(L^{G})$ corresponding to the largest $r$ eigenvalues are the left-most $r$ vectors in $V$.
For the sake of convenience, in this paper, we simply refer to such filter as a low pass graph filter. A low pass graph signal is defined as a graph signal generated by the excitation of a low pass filter, including but not limited to white noise. Using \Cref{deflowpass}, the ideal low pass graph filter has $\eta_{r}=0$. Low pass graph filters are commonly employed in modeling practical network data, as detailed in \cite{lowpass}. Below we present two examples.


\begin{example}{\rm\cite{diffusion, IIR}}\label{ep1}
The graph filter $\mathcal{H}(L^{G})$ defined by 
$$\mathcal{H}(L^{G})=\sum\limits_{t=0}^{+\infty}\frac{\sigma^{t}}{t!}(L^{G})^{t}$$ is a strong low pass graph filter for large $\sigma$. Corresponding to the above \eqref{H(L)}, we know that $a_{t}=\frac{\sigma^{t}}{t!}$ and $T=+\infty$. According to the definition of the series, we can also write it as $\mathcal{H}(L^{G})=e^{-\sigma L^{G}}$. Note that 
$$\eta_{r}=\frac{e^{-\sigma\lambda_{r+1}}}{e^{-\sigma\lambda_{r}}}=e^{-\sigma(\lambda_{r+1}-\lambda_{r})}.$$
Since $\lambda_{r+1}>\lambda_{r}$ and $\sigma>0$, the low pass ratio decreases to zero exponentially with $\sigma$.
\end{example}

\begin{example}{\rm\cite{IIR}}\label{ep2}
The graph filter $\mathcal{H}(L^{G})$ defined by 
$$\mathcal{H}(L^{G})=\sum\limits_{t=0}^{+\infty}(-\alpha)^{t}(L^{G})^{t}$$ with $\alpha>0$, is a weak low pass graph filter for small $\alpha$. Corresponding to the above \eqref{H(L)}, we know that $a_{t}=(-\alpha)^{t}$ and $T=+\infty$. We can also express the above graph filter as $\mathcal{H}(L^{G})=(I+\alpha L^{G})^{-1}$. Note that $$\eta_{r}=\frac{1+\alpha\lambda_{r}}{1+\alpha\lambda_{r+1}}=1-\alpha\frac{\lambda_{r+1}-\lambda_{r}}{1+\alpha\lambda_{r}}.$$
Since $\lambda_{r+1}>\lambda_{r}$ and $\alpha>0$, the low pass ratio $\eta_{r}\approx1$ when $\alpha$ is sufficiently small.
\end{example}

\section{Blind Extraction with Low Pass Graph Filters}\label{sec3}

In this section, we establish the interesting connection between EEP and low pass graph signals. Moreover, we propose \Cref{alg:2}, namely BE-EEPs, to address this blind extraction problem by establishing the matrix optimization model for extraction of EEPs via low pass graph filters.


\subsection{Problem Statement}\label{sec3-1}
Given the observations $\textbf{y}^{1}, \cdots, \textbf{y}^{m}\in\mathbb{R}^{n}$, which are the graph filter's outputs subject to the unknown excitation $\textbf{x}^{1}, \cdots, \textbf{x}^{m}\in\mathbb{R}^{n}$
\begin{equation}\label{1}
\textbf{y}^{l}=\mathcal{H}(L^{G})\textbf{x}^{l}+\omega^{l},\ l=1, \cdots, m,
\end{equation}
where $\mathcal{H}(L^{G})$ is defined by \eqref{H(L)}, $\omega^{l}\in\mathbb{R}^{n}$ represents the noise satisfying normal distribution $\omega^{l}\sim\mathcal{N}(0, \sigma^{2}I)$, our aim is to find an EEP $\pi$ of $G$ from the graph filters. The following assumptions are standard assumptions based on which will be convenient for our analysis.

\begin{assumption}\label{assum1}
The graph filter $\mathcal{H}(L^{G})$ of order $T$ is a low pass graph filter.
\end{assumption}

\begin{assumption}\label{assum2}\footnote{Most prior work on graph learning made this assumption about the excitation graph signals $\textbf{x}^{l}$, $l=1,\cdots,m$, such as in \cite{no22, no25, no26, no32}. This assumption implies that the excitations at each node are independent.}
The excitation graph signals $\textbf{x}^{l}$ is zero-mean white noise with $\mathbb{E}[\textbf{x}^{l}(\textbf{x}^{l})^{\top}]=I$, $l=1, \cdots, m$.
\end{assumption}


With \Cref{assum1} and \Cref{assum2}, the covariance matrix of the observed signals $\textbf{y}^{l}$, is given by
\begin{equation}\label{3-1}
\Sigma:=\mathbb{E}[\textbf{y}^{l}(\textbf{y}^{l})^{\top}]=\mathcal{H}(L^{G})\mathcal{H}(L^{G})^{\top}+\sigma^{2}I.
\end{equation}

To derive EEP of $G$ from low pass graph filters, the indicator matrix $H\in\mathbb{R}^{n\times r}$ can be used to derive the EEP $\pi$. Take the graph in \Cref{fig-quotient} (a) as an example. We clearly know that
$$
L^{G}=\left[\begin{array}{c c : c c c : c c c c c c }
3&0&-1&-1&-1&0&0&0&0&0&0\\
0&3&-1&-1&-1&0&0&0&0&0&0\\
\hdashline
-1&-1&4&0&0&-1&-1&0&0&0&0\\
-1&-1&0&5&-1&0&0&-1&-1&0&0\\
-1&-1&0&-1&5&0&0&0&0&-1&-1\\
\hdashline
0&0&-1&0&0&2&-1&0&0&0&0\\
0&0&-1&0&0&-1&2&0&0&0&0\\
0&0&0&-1&0&0&0&1&0&0&0\\
0&0&0&-1&0&0&0&0&1&0&0\\
0&0&0&0&-1&0&0&0&0&1&0\\
0&0&0&0&-1&0&0&0&0&0&1
\end{array}\right],\
H=\left[\begin{array}{c c c}
1&0&0\\
1&0&0\\
\hdashline
0&1&0\\
0&1&0\\
0&1&0\\
\hdashline
0&0&1\\
0&0&1\\
0&0&1\\
0&0&1\\
0&0&1\\
0&0&1
\end{array}\right],
$$
and
$$
L^{G/\pi}=\left[\begin{array}{c c c}
3&-3&0\\
-2&4&-2\\
0&-1&1
\end{array}\right].
$$
One can verify that $L^{G}H=HL^{G/\pi}$. Then we can calculate that $\mathbb{C}=\{C_{1}, C_{2}, C_{3}\}$ and $C_{1}=\{1,2\}$, $C_{2}=\{3,4,5\}$, $C_{3}=\{6,7,8,9,10,11\}$.
The details are summarized in \Cref{alg:1}.
\begin{algorithm}[H]
\caption{The procedure to finding EEP $\pi$ through an indicator matrix $H$}
\label{alg:1}
\begin{algorithmic}[1]
\REQUIRE an indicator matrix $H\in\mathbb{R}^{n\times r}$, an $n$-dimensional zero vector ${\bf c}=\{c_{1},\cdots,c_{n}\}$.
\STATE For $i=1,2,\cdots,n$, set $c_{i}=k$ if $H_{ik}\neq0$, $k=1,2,\cdots,r$.
\STATE $C_{k}=\{i:\ c_{i}=k\}$ where $k=1,2,\cdots,r$ and $i=1,2,\cdots,n$.
\ENSURE an EEP $\pi$ with $\mathbb{C}=\{C_{1}, C_{2}, \cdots, C_{r}\}$.
\end{algorithmic}
\end{algorithm}

\subsection{Connections between EEP and low pass graph filters}\label{sec3-2}
Note that as mentioned in \Cref{sec3-1}, EEP of a graph can be derived by the indicator matrix $H$. In this part, we will derive the relation between the indicator matrix $H$ and the low pass graph filters. The rough idea is as follows. We first establish the relationship between $\mathcal{H}(L^{G})$ and $H$ in \Cref{remark-1}. Combining \Cref{prop-ep-eep}, we derive that $\mathcal{H}(L^{G})$ also possesses structural eigenvectors as described in \Cref{prop-ep-eep}(ii). Furthermore, through \Cref{prop1}, we infer that starting from the covariance matrix, we can obtain the indicator matrix $H$ and consequently derive the EEP $\pi$ with $\mathbb{C}=\{C_{1}, C_{2}, \cdots, C_{r}\}$.

\begin{proposition}\label{remark-1}
For an EEP $\pi$ as indicated by $H$, it holds that $\mathcal{H}(L^{G})H=H\mathcal{H}(L^{G/\pi})$.
\end{proposition}
{\bf Proof.}
Note that $\mathcal{H}(L^{G})$ is a polynomial of the matrix $L^{G}$. According to $L^{G}H=HL^{G/\pi}$ in \Cref{prop-ep-eep}(i), it holds that $(L^{G})^{r}H=H(L^{G/\pi})^{r}$ for any nonnegative integer $r$. Therefore, by \eqref{H(L)}, it holds that $$\mathcal{H}(L^{G})H=\sum\limits^{T}_{t=0}a_{t}(L^{G})^{t}H=\sum\limits^{T}_{t=0}a_{t}H(L^{G/\pi})^{t}=H\mathcal{H}(L^{G/\pi}).$$
\hfill$\Box$

Since the eigenvectors of $L^{G}$ are also eigenvectors of the polynomial of the matrix $L^{G}$, that is, the structural eigenvectors of $L^{G}$ are also the eigenvectors of $\mathcal{H}(L^{G})$.
\begin{definition}
The structural eigenvectors of $L^{G}$ are defined as the structural eigenvectors of $\mathcal{H}(L^{G})$.
\end{definition}

The following result reveals the equality relationship between the structural eigenvectors of $\mathcal{H}(L^{G})$ and some of the eigenvectors of the covariance matrix $\Sigma$.

\begin{proposition}\label{prop1}
For an EEP $\pi$ as indicated by $H$, the eigenvectors of $\Sigma$ in the form of $H\tilde{v}$ are the structural eigenvectors of $\mathcal{H}(L^{G})$, where $\tilde{v}$ is an eigenvector of $\mathcal{H}(L^{G/\pi})$.
\end{proposition}
{\bf Proof.}
Let $\mathcal{H}(L^{G/\pi})\tilde{v}=\tilde{\lambda} \tilde{v}$, where $\tilde{v}$ is an eigenvector of $\mathcal{H}(L^{G/\pi})$ with eigenvalue $\tilde{\lambda}$, we obtain that
\begin{equation}
\begin{aligned}
\label{hv-stronglowpass}
\Sigma H\tilde{v}&= \mathcal{H}(L^{G})(\mathcal{H}(L^{G}))^{\top}H\tilde{v}+\sigma^{2}I H\tilde{v} \\
&= \mathcal{H}(L^{G})\mathcal{H}(L^{G})H\tilde{v}+\sigma^{2}I H\tilde{v}\\
&= \mathcal{H}(L^{G})H\mathcal{H}(L^{G/\pi})\tilde{v}+\sigma^{2}I H\tilde{v} \\
&= H(\mathcal{H}(L^{G/\pi}))^{2}\tilde{v}+\sigma^{2}I H\tilde{v} \\
&= (\tilde{\lambda}^{2}+\sigma^{2})H\tilde{v},
\end{aligned}
\end{equation}
where the second equality stems from the fact that $\mathcal{H}(L^{G})$ is a symmetric matrix; the third equality and the forth equality are from \Cref{remark-1}; the fifth equality comes from $\mathcal{H}(L^{G/\pi})\tilde{v}=\tilde{\lambda} \tilde{v}$. Hence, the eigenvectors of $\Sigma$ in the form of $H\tilde{v}$ are the structural eigenvectors of $\mathcal{H}(L^{G})$.
\hfill$\Box$

\begin{definition}
Let $v$ be an eigenvector of $\Sigma$. For an EEP $\pi$ as indicated by $H$, if there exists an eigenvector of $\mathcal{H}(L^{G/\pi})$, denoted as $\tilde{v}$, such that $v=H\tilde{v}$, then $v$ is referred to as a structural eigenvector of $\Sigma$.
\end{definition}

Through the two propositions above, for low pass graph signal models that satisfy \Cref{assum1} and \Cref{assum2}, the covariance matrix contains structural eigenvectors that reflect the indicator matrix. Therefore, by computing the eigenvectors of the covariance matrix and identifying the structural eigenvectors among them, we can determine the indicator matrix $H$ and subsequently find the EEP. The process of identifying these structural eigenvectors is described in \Cref{sec3.4}.
\subsection{A Function Related to Structural Eigenvectors and an Indicator Matrix}\label{sec3-3}
For convenience, we perform the spectral decomposition of $\mathcal{H}(L^{G})$ and $\Sigma$ respectively. Let $\overline{V}\in\mathbb{R}^{n\times r}$ be a matrix containing $r$ structural eigenvectors of $\mathcal{H}(L^{G})$. Let $P\in\mathbb{R}^{n\times r}$ be a matrix containing $r$ structural eigenvectors of $\Sigma$, and is associated with $r$ eigenvalues $\xi_{1}\geq\cdots\geq\xi_{r}$. Then we can get that $\overline{V}=P$ due to \Cref{remark-1} and \Cref{prop1}.
Define the cost function \cite{no26} as follows
\begin{equation}
\label{F}
F(\mathbb{C}, V):=\sum\limits_{k=1}^{r}\sum\limits_{i\in C_{k}}\|V_{i,\cdot}-\frac{1}{|C_{k}|}\sum\limits_{j\in C_{k}}V_{j,\cdot}\|^{2}_{2},
\end{equation}
where $V_{i,\cdot}$ denotes the $i$-th row of $V$. For an EEP, the components of any structural eigenvector corresponding to the same class do not change. Therefore, for an EEP $\overline{\pi}$ of $G$ with $\overline{\mathbb{C}}=\{\overline{C}_{1}, \overline{C}_{2}, \cdots, \overline{C}_{r}\}$, we obtain that $F(\overline{\mathbb{C}}, \overline{V})=F(\overline{\mathbb{C}}, P)=0$.

To reformulate \eqref{F} in the form of a binary indicator matrix $H\in\mathbb{R}^{n\times r}$ of the partition $\pi$ with $\mathbb{C}$, we define a normalized indicator matrix (see \cite{nim})
\begin{equation}\label{hatH}
\widehat{H}=H{\rm Diag}\left(\frac{1}{\sqrt{|C_{1}|}}, \frac{1}{\sqrt{|C_{2}|}}, \cdots, \frac{1}{\sqrt{|C_{r}|}}\right),
\end{equation}
where $|C_{i}|$ denote the number of vertices in the set $C_{i}$. It holds that
\begin{equation}
\left(\widehat{H}\widehat{H}^{\top}\right)_{ij} =\left\{\begin{array}{ll}
		\frac{1}{|C_{k}|},&{\rm if}\ i\in C_{k} \ {\rm and}\ j\in C_{k}, k=1,2,\cdots,r,\\
		0,&{\rm otherwise}.
		\end{array}	\right.
\end{equation}
It gives the following
\begin{equation}
\label{F1}
\sum\limits_{k=1}^{r}\sum\limits_{i\in C_{k}}\|V_{i,\cdot}-\frac{1}{|C_{k}|}\sum\limits_{j\in C_{k}}V_{j,\cdot}\|^{2}_{2}=\|V-\widehat{H}\widehat{H}^{\top}V\|^{2}_{F}.
\end{equation}
Let $\overline{H}$ be the normalized indicator matrix corresponding to the EEP $\overline{\pi}$. Then we know that $\|\overline{V} - \overline{H}(\overline{H})^{\top}\overline{V}\|_{F}^{2} = 0$.
Below, we will use \eqref{F1} to relate the structural eigenvectors and an indicator matrix to construct an optimization model.
\subsection{The Optimization Model}\label{sec3.4}
As mentioned above, our goal is to derive the EEP of the graph $G$ from the outputs of an unknown network process represented by a low pass graph filter.
Given the observed signals $\textbf{y}^{l}$, $l=1,2,\cdots,m$, instead of the covariance matrix $\Sigma$, the following sample covariance matrix (denoted by $\widehat{\Sigma}$)

\begin{equation}\label{hatSigma}
\widehat{\Sigma}=\frac{1}{m}\sum\limits_{l=1}^{m}\textbf{y}^{l}(\textbf{y}^{l})^{\top}
\end{equation}
is available in practice.

According to Proposition 1 in \cite{no26} and \Cref{assum1}, most of top $r$ eigenvectors of $\widehat{\Sigma}$ are structural eigenvectors. Hence we can simply choose the top $r$ eigenvectors of $\widehat{\Sigma}$ as the structural eigenvectors. Subsequently, we will utilize these top $r$ eigenvectors of $\widehat{\Sigma}$ to approximately extract the EEP.

Let $\widehat{P}\in\mathbb{R}^{n\times r}$ be a matrix containing $r$ estimated structural eigenvectors of $\widehat{\Sigma}$ with $r$ eigenvalues $\hat{\xi}_{1}\geq\cdots\geq\hat{\xi}_{r}$.
Recall the cost function in \eqref{F}. We minimize the cost function with respect to the estimated structural eigenvectors $\widehat{P}$, which gives the following model:
\begin{equation}
\label{F-hat}
\min\limits_{\mathbb{C}}f(\mathbb{C}, \widehat{P}):=\min\limits_{\mathbb{C}}\sum\limits_{k=1}^{r}\sum\limits_{i\in C_{k}}\|\widehat{P}_{i,\cdot}-\frac{1}{|C_{k}|}\sum\limits_{j\in C_{k}}\widehat{P}_{j,\cdot}\|^{2}_{2}.
\end{equation}
Keeping \eqref{F1} in mind, for the graph signal model with low pass graph filters, we can search for an approximate EEP by solving the following optimization problem (recall that $\widehat{P}$ is given)
\begin{equation}
\label{F2}
\min\limits_{\widehat{H}\in\mathbb{H}_{r}}\|\widehat{P}-\widehat{H}\widehat{H}^{\top}\widehat{P}\|^{2}_{F},
\end{equation}
where $$\mathbb{H}_{r}=\{\widehat{H}\in\mathbb{R}^{n\times r} \ | \ \widehat{H}^{\top}\widehat{H}=I_{r},\ \widehat{H}(\widehat{H})^{\top}{\bf 1}_{n}={\bf 1}_{n},\ \|\widehat{H}_{i,:}\|_{0}=1,\ \text{{\rm and }} \ \widehat{H}\geq0\}.$$ The constraint $\|\widehat{H}_{i,:}\|_{0}=1$ shows that $\widehat{H}$ have exactly one non-zero element per row, which denotes the clustering membership. The constraints $\widehat{H}\geq0$, $\widehat{H}^{\top}\widehat{H}=I_{r}$ and $\widehat{H}(\widehat{H})^{\top}{\bf 1}_{n}={\bf 1}_{n}$ indicate that $\widehat{H}$ is a normalized indicator matrix. $\widehat{H}^{\top}\widehat{H}=I_{r}$ ensures that each column of $\widehat{H}$ contains non-zero elements. Furthermore, $\widehat{H}(\widehat{H})^{\top}{\bf 1}_{n}={\bf 1}_{n}$ ensures that the non-zero elements in the columns of $\widehat{H}$ are consistent.
The following algorithm, named as BE-EEPs, summarizes the above procedure in a clear framework.
\begin{algorithm}[H]
\caption{BE-EEPs: Blind EEPs Extraction by low pass graph signal}
\label{alg:2}
\begin{algorithmic}[1]
\REQUIRE Graph signals $\{\textbf{y}^{l}\}^{m}_{l=1}$; desired number of classes $r$.
\STATE Compute the sample covariance $\widehat{\Sigma}$ as \eqref{hatSigma}.
\STATE Find the top $r$ eigenvectors of $\widehat{\Sigma}$, denoted as $\widehat{P}\in\mathbb{R}^{n\times r}$.
\STATE Solve \eqref{F2} to get $\widehat{H}$.
\STATE Use \Cref{alg:1} to get $\widehat{\mathbb{C}}$.
\ENSURE $\widehat{H}$ and $\widehat{\mathbb{C}}=\{\widehat{C}_{1}, \widehat{C}_{2}, \cdots, \widehat{C}_{r}\}$.
\end{algorithmic}
\end{algorithm}

\section{Error Bound Analysis}\label{sec4}
In this part, we derive theoretical error bounds for the performance of \Cref{alg:2} under certain assumptions. We need the following result, which comes from Theorem 2 in \cite{no26}.

\begin{theorem}\label{thm-A1}{\rm \cite[Theorem 2]{no26}}
Let $\textbf{y}^{1}$, $\cdots$, $\textbf{y}^{m}$ be independent samples of the graph filter as in \eqref{1}. Assume that for some $K\geq1$, $\|\textbf{y}^{l}\|_{2}\leq K(\mathbb{E}(\|\textbf{y}^{l}\|_{2}^{2}))^{\frac{1}{2}}$, $l=1,\cdots,m$ almost surely. Let the effective rank of the covariance be $b=\frac{tr(\Sigma)}{\|\Sigma\|_{2}}$, then for every $c'>0$ it holds that:
\begin{equation}
\|\widehat{\Sigma}-\Sigma\|_{2}\leq \epsilon\left(\sqrt{\frac{K^{2}b{\rm log}(\frac{n}{c'})}{m}}+\frac{K^{2}b{\rm log}(\frac{n}{c'})}{m}\right)\|\Sigma\|_{2}
\end{equation}
with probability at least $1-c'$, where $\epsilon$ is an constant.
\end{theorem}

\Cref{thm-A1} bounds the difference between sample covariance and covariance. Following this, and based on \Cref{thm-A1}, we provide the upper bound for the difference between the cost function $F(\overline{\mathbb{C}}, \overline{V})$ under ideal conditions and the cost achieved by the partition $\widehat{\mathbb{C}}=\{\widehat{C}_{1}, \widehat{C}_{2}, \cdots, \widehat{C}_{r}\}$ identified by \Cref{alg:2}, as demonstrated in \Cref{thm-A2}.

\begin{theorem}\label{thm-A2}
Let $\textbf{y}^{1}$, $\cdots$, $\textbf{y}^{m}$ be independent samples of the graph filter as in \eqref{1}. Let $b=\frac{tr(\Sigma)}{\|\Sigma\|_{2}}$ and the following conditions hold:

{\rm (\romannumeral1)} \Cref{alg:2} finds a solution $\widehat{\mathbb{C}}=\{\widehat{C}_{1}, \widehat{C}_{2}, \cdots, \widehat{C}_{r}\}$ that exactly minimizes the problem \eqref{F2};

{\rm (\romannumeral2)} There exists $\delta>0$ such that $\|\widehat{\Sigma}-\Sigma\|_{2}+\delta\leq \xi_{r}-\hat{\xi}_{r+1}$;

{\rm (\romannumeral3)} For some $K\geq1$, $\|\textbf{y}^{l}\|_{2}\leq K(\mathbb{E}(\|\textbf{y}^{l}\|_{2}^{2}))^{\frac{1}{2}}$, $l=1,\cdots,m$ almost surely.

Then for every $c'>0$, it holds that
\begin{equation}\label{thm3.2-1}
\sqrt{F(\widehat{\mathbb{C}}, \overline{V})}-\sqrt{F(\overline{\mathbb{C}}, \overline{V})}\leq2\sqrt{2r}\frac{\epsilon\left(\sqrt{\frac{K^{2}b{\rm log}(\frac{n}{c'})}{m}}+\frac{K^{2}b{\rm log}(\frac{n}{c'})}{m}\right)\|\Sigma\|_{2}}{\delta}
\end{equation}
with probability at least $1-c'$, where $\epsilon$ is an constant.
\end{theorem}
{\bf Proof.}
It follows from \eqref{F1} and $F(\overline{\mathbb{C}}, \overline{V})=0$ that
\begin{equation*}
\sqrt{F(\widehat{\mathbb{C}}, \overline{V})}-\sqrt{F(\overline{\mathbb{C}}, \overline{V})}=\|\overline{V}-\widehat{H}\widehat{H}^{\top}\overline{V}\|_{F}.
\end{equation*}
The rest of the proof is similar to that in the proof of Theorem $1$ in \cite{no26}.
\hfill$\Box$
\eqref{thm3.2-1} implies that the projection error associated with the estimated structural eigenvectors $\widehat{P}$ will be small if the sample size is large enough.

\section{Three Methods for Solving \eqref{F2}}\label{sec5}
In this part, we will discuss how to solve \eqref{F2} in \Cref{alg:2} using three methods.

To be precise, \eqref{F2} can also be written in the following matrix optimization problem (recall that $\widehat{P}$ is given)
\begin{equation}
\label{F2-1}
\begin{aligned}
\min\limits_{\widehat{H}\in\mathbb{R}^{n\times r}}& \hspace{2mm}\|\widehat{P}-\widehat{H}\widehat{H}^{\top}\widehat{P}\|^{2}_F\\
{\rm s.t.}&\hspace{2mm}\widehat{H}^{\top}\widehat{H}=I_{r},\\
&\hspace{2mm}\widehat{H}(\widehat{H})^{\top}{\bf 1}_{n}={\bf 1}_{n},\\
&\hspace{2mm}\|\widehat{H}_{i,:}\|_{0}=1,\\
&\hspace{2mm}\widehat{H}\geq0.
\end{aligned}
\end{equation}
The objective function in \eqref{F2-1} is continuously differentiable.
We know that the feasible set of \eqref{F2-1} is non-convex with some combinatorial structures. 
Due to the combinatorial structures, solving \eqref{F2-1} is generally NP-hard. This model also arises in various fields such as classification \cite{classifi1, onmf}, gene expression studies \cite{Gene}, and text mining \cite{text}. We summarize two key points of \eqref{F2-1} that need special attention: (\romannumeral1) $\widehat{P}$ may have negative elements; (\romannumeral2) each row of $\widehat{H}$ must contain exactly one positive element. While optimization with nonnegative orthogonal constraints has been extensively explored \cite{opnmf, no33, onmf-1, onmf-2, onmf-3}, few studies consider models that simultaneously satisfy these two points. Taking these into account, we choose to address the problem using the following three methods.

\paragraph{Method One} While the classical K-means model is typically represented using $k$ centroids, it can also be expressed using an indicator matrix. As is well known, the classical K-means model can also be rewritten as \eqref{F2-1} \cite{nim, k-means2}.
It is evident that to solve problem $\eqref{F2-1}$, we only need to solve the K-means problem. Additionally, \cite{2-4eep} demonstrated a significant duality between the extraction problem of the EEP and the K-means problem on the rows of Laplacian eigenvectors. Therefore, we can consider using the K-means algorithm to solve problem $\eqref{F2-1}$.
There are many software packages available to solve the K-means model, such as the built-in K-means clustering function in MATLAB\footnote{Details can be found at https://ww2.mathworks.cn/help/stats/kmeans.html} and `{\bf sklearn.cluster.KMeans}' in Python\footnote{Details can be found at https://scikit-learn.org/stable/modules/clustering.html\#k-means}. Here, we choose to use the `{\bf sklearn.cluster.KMeans}' class from the popular Python machine learning library scikit-learn to solve the K-means model.
\paragraph{Method Two} On the other hand, the orthonormal projective nonnegative matrix factorization (OPNMF) model in \cite{opnmf} takes the following form 
\begin{equation}
\label{opnmf}
\begin{aligned}
\min\limits_{\widehat{H}\in\mathbb{R}^{n\times r}}& \hspace{2mm}\|\bar{P}-\widehat{H}\widehat{H}^{\top}\bar{P}\|^{2}_{F}\\
{\rm s.t.}&\hspace{2mm}\widehat{H}^{\top}\widehat{H}=I_{r},\\
&\hspace{2mm}\widehat{H}\geq0,
\end{aligned}
\end{equation}
where $\bar{P}\in\mathbb{R}^{n\times r}_{+}$ is given. \eqref{opnmf} is equivalent to the classical orthonormal nonnegative matrix factorization (ONMF) model \cite{onmf} and has been widely studied, such as in \cite{onmf,opnmf,pnmf-1,pnmf-2,pnmf-3}. Many approaches have been proposed to solve ONMF or OPNMF \cite{opnmf,pnmf-4,pnmf-5}. 
Comparing \eqref{opnmf} with \eqref{F2-1}, the difference lies in two aspects: One is that the constraints $\|\widehat{H}_{i,:}\|_{0}=1$ and $\widehat{H}(\widehat{H})^{\top}{\bf 1}_{n}={\bf 1}_{n}$ are dropped in \eqref{opnmf}; The other is $\bar{P}$ in \eqref{opnmf} is nonnegative whereas $\widehat{P}$ in \eqref{F2-1} may have negative elements. By letting $\bar{P}=\widehat{P}-a{\bf1}_{n\times r}$, where $a$ is the minimum value in the matrix $\widehat{P}$, one can guarantee that $\bar{P}$ is nonnegative. Note that, due to the existence of \eqref{hv-stronglowpass} and \Cref{thm-A1}, and considering that the eigenvectors corresponding to different eigenvalues of the Laplacian matrix are orthogonal, $a$ must necessarily be a negative value. Furthermore, due to the extraction concept of EEP in \Cref{sec3-3}, this transformation will not affect the original objective function value. Next, we choose the practical exact penalty method EP4Orth+ proposed in \cite{no33} since it performs well and can always return high quality orthogonal nonnegative matrices, as demonstrated in \cite{no33}. Additionally, the postprocessing step of EP4Orth+ ensures that the resulting matrix satisfies the constraint $\|\widehat{H}_{i,:}\|_{0}=1$. We refer to \cite{no33} for more details of the method.

\paragraph{Method Three} Noting that $\widehat{P}$ in \eqref{F2-1} may have negative elements, inspired by \cite{opnmf}, we choose to rewrite the version without the zero norm constraint of \eqref{F2-1} as a projective semi-nonnegative matrix factorization (Projective Semi-NMF) problem. We then solve it using an iterative Lagrangian approach \cite{opnmf}, detailed as follows. We first remove the zero norm constraint and $\widehat{H}\widehat{H}^{\top}{\bf 1}_{n}={\bf 1}_{n}$ from \eqref{F2-1}, resulting in the following form
\begin{equation}
\label{F2-nozero}
\begin{aligned}
\min\limits_{\widehat{H}\in\mathbb{R}^{n\times r}}& \hspace{2mm}\|\widehat{P}-\widehat{H}\widehat{H}^{\top}\widehat{P}\|^{2}_F\\
{\rm s.t.}&\hspace{2mm}\widehat{H}^{\top}\widehat{H}=I_{r},\\
&\hspace{2mm}\widehat{H}\geq0.
\end{aligned}
\end{equation}

Let $\widehat{p}_{i}:=\widehat{P}_{i,\cdot}\in \mathbb{R}^{r}$, $i=1,\cdots,n$. As shown in sections II-E and II-F of \cite{opnmf}, let
$$\Phi=\left[\phi^{\top}(\widehat{p}_{1}), \phi^{\top}(\widehat{p}_{2}),  \cdots, \phi^{\top}(\widehat{p}_{n})\right]^{\top},$$
where $\phi$ is a vector function. Denote the kernel matrix \cite{opnmf} $K=\Phi\Phi^{\top}$, with $K_{ij}=\phi(\widehat{p}_{i})\phi^{\top}(\widehat{p}_{j})$, $i,j=1,\cdots,n$. Here we choose the linear kernel function $\phi(\widehat{p}_{i}) = \widehat{p}_{i}$ for $i=1,\cdots,n$, i.e., $\Phi=\widehat{P}$. Hence, we can know objective function $\|\widehat{P}-\widehat{H}\widehat{H}^{\top}\widehat{P}\|^{2}_F$ in \eqref{F2-nozero} is equal to $\|\Phi-\widehat{H}\widehat{H}^{\top}\Phi\|^{2}_F$. Because $\widehat{P}$ may not be nonnegative, $K$ is not necessarily nonnegative either. As shown in \cite{opnmf}, for a kernel matrix $K$ that has both positive and negative elements, we can separate the positive and negative parts of the matrix by calculating
$$K_{ij}^{+}=\frac{|K_{ij}|+K_{ij}}{2} \ {\rm and}\ K_{ij}^{-}=\frac{|K_{ij}|-K_{ij}}{2},\ i,j=1,\cdots,n,$$
where $K=K^{+}-K^{-}$ with $K^{+}\geq0$ and $K^{-}\geq0$. Next, we rewrite \eqref{F2-nozero} as a projective Semi-NMF problem 
\begin{equation}
\label{psnmf}
\begin{aligned}
\min\limits_{\widehat{H},\widehat{U}\in\mathbb{R}^{n\times r}}& \hspace{2mm}\|\Phi-\widehat{H}\widehat{U}^{\top}\Phi\|^{2}_F\\
{\rm s.t.}&\hspace{2mm}\widehat{U}=\widehat{H},\\
&\hspace{2mm}\widehat{H}^{\top}\widehat{H}=I_{r},\\
&\hspace{2mm}\widehat{H}\geq0.
\end{aligned}
\end{equation}
Then, following the derivation in Appendix VIII of \cite{opnmf}, we obtain its iterative Lagrangian solution
\begin{equation}
\label{iterative}
\widehat{H}^{(k+1)}_{ij}=\widehat{H}^{(k)}_{ij}\frac{\left(K^{+}\widehat{H}^{(k)}+\widehat{H}^{(k)}(\widehat{H}^{(k)})^{\top}K^{-}\widehat{H}^{(k)}\right)_{ij}}{\left(K^{-}\widehat{H}^{(k)}+\widehat{H}^{(k)}(\widehat{H}^{(k)})^{\top}K^{+}\widehat{H}^{(k)}\right)_{ij}}.
\end{equation}
It is proved in \cite{opnmf} that the iterative scheme in \eqref{iterative} is an iterative Lagrangian solution of \eqref{psnmf}.

However, unfortunately, while this method is suitable for cases where $\widehat{P}$ contains negative elements, it has certain limitations. Specifically, the $\widehat{H}$ obtained through the above method does not satisfy the zero norm constraint $\|\widehat{H}_{i,:}\|_{0}=1$. Generally, the following two situations may arise: (\romannumeral1) there may be more than one nonzero element in a row of $\widehat{H}$; (\romannumeral2) iteration \eqref{iterative} cannot restore zero entries to positive values, which may result in rows in $\widehat{H}$ being all zeros. This also leads to the fact that the iterative \eqref{iterative} does not guarantee orthogonality. To address these issues, we implement some corrective measures.

For case (\romannumeral1), inspired by \cite{no33}, we extract the index with the maximum value in a row. As for case (\romannumeral2), we need to stop the iteration immediately upon discovering a row in $\widehat{H}$ with all zero elements. Then, we apply the operation for (\romannumeral1) to ensure that $\widehat{H}$ satisfies the zero norm constraint in \eqref{F2-1}. Notably, for the blind extraction problem of EEP, our focus lies in identifying the correct EEP rather than the specific values of $\widehat{H}$. This motivation also drives us to perform the aforementioned operations.

\section{Numerical Results}\label{sec6}

In this section, we show the numerical performance of \Cref{alg:2} to verify the efficiency of our proposed approach. As mentioned in \Cref{sec5}, we compare the K-means algorithm, the practical exact penalty method and the iterative Lagrangian approach, denoted as K-means, EP4Orth+ and PSNMF, to solve \eqref{F2-1}. The experiments are performed in Windows 10 on an Intel(R) Core(TM) i7-1065G7 CPU at 1.30 GHZ with 16.0 GB of RAM. All codes are written in MATLAB R2022a.

We evaluate the performances of algorithms using synthetic graph signals. The graph sampled from a symmetric hierarchical network model \cite{2-4eep} are considered, which guarantee that the sampled graph has an EEP. For more details of this model we refer the reader to \cite{2-4eep}. Here, we set the graph to have $378$ nodes and an EEP with three same sized classes. Furthermore, the observed graph signal is generated as \eqref{1}. Let $D_{max}$ be the highest degree of graph $G$. We consider the two types of low pass graph filters mentioned in \Cref{ep1} and \Cref{ep2}: (a) $\mathcal{H}_{strong}(L^{G})=e^{-\sigma L^{G}}$ where $\sigma=10/D_{max}$; (b) $\mathcal{H}_{weak}(L^{G})=(I+\alpha L^{G})^{-1}$ where $\alpha=0.5/D_{max}$. Notice
that $\mathcal{H}_{strong}(L^{G})$ is strong low pass with $\eta_{r}\approx0$, and $\mathcal{H}_{weak}(L^{G})$ is weak low pass with $\eta_{r}\approx1$.

\begin{figure*}[ht]
\begin{center}
\subfigure{
\includegraphics[width=2.6in,height=2in]{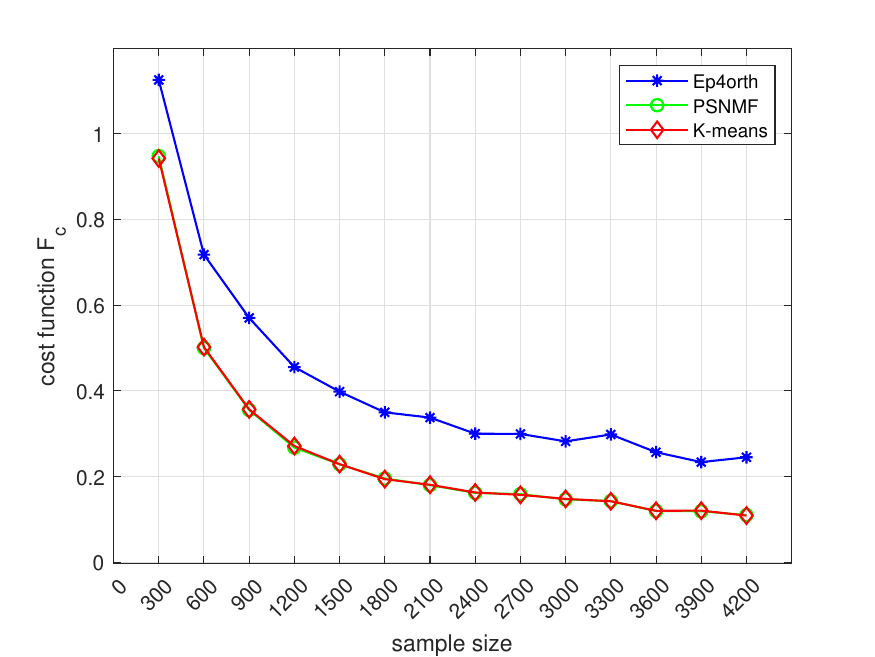}
}
\subfigure{
\includegraphics[width=2.6in,height=2in]{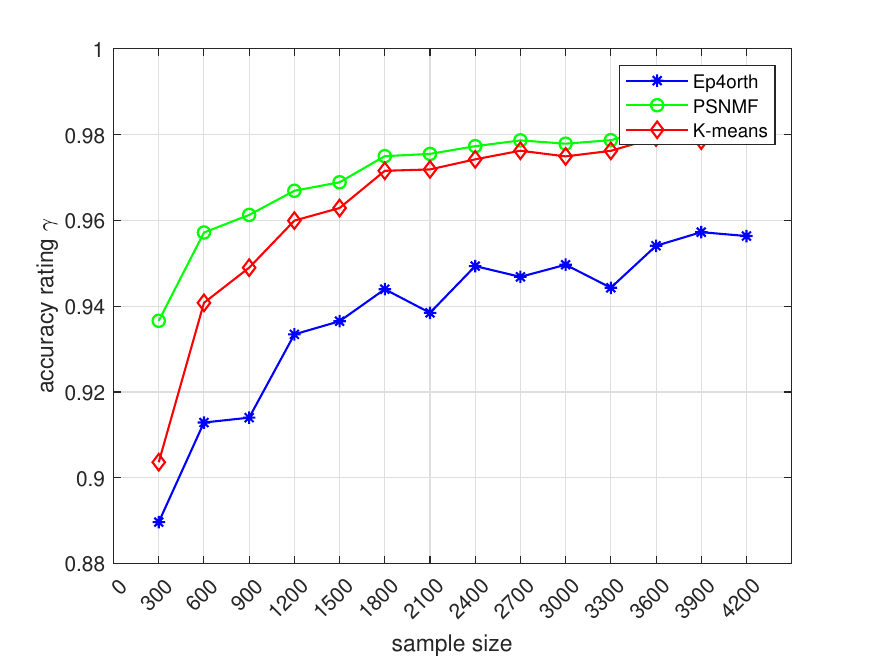}
}
\caption{The cost function $F_{c}$ and the group accuracy rate $\gamma$ of K-means, PSNMF and EP4Orth+ for strong low pass graph signals.}
\label{fig1}
\end{center}
\end{figure*}

\begin{figure*}[ht]
\begin{center}
\subfigure{
\includegraphics[width=2.6in,height=2in]{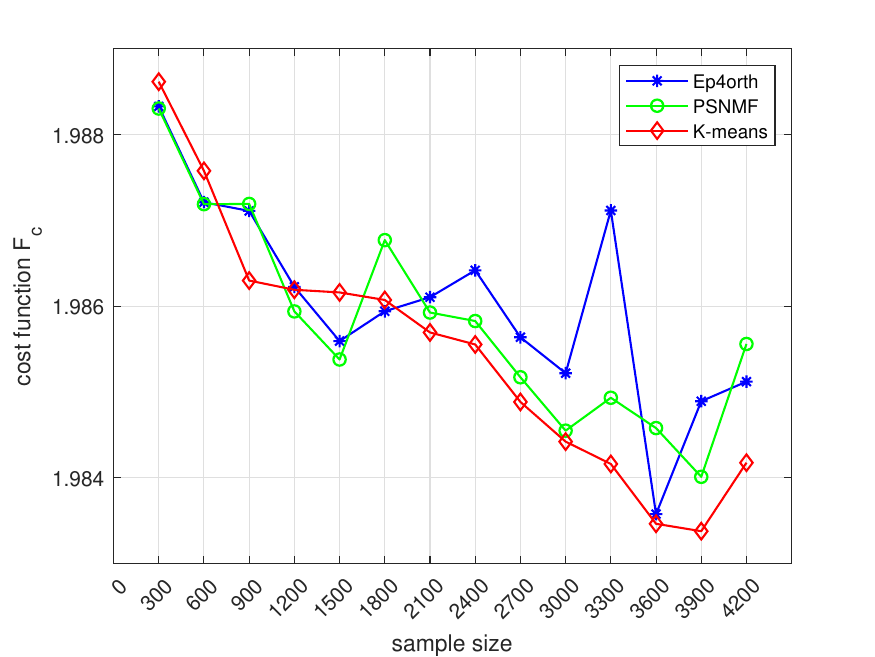}
}
\subfigure{
\includegraphics[width=2.6in,height=2in]{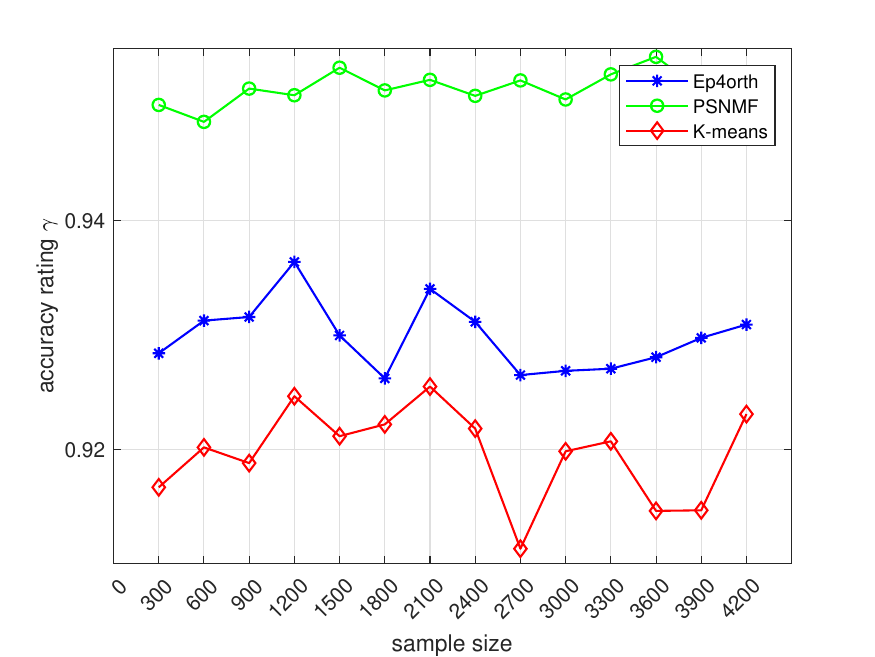}
}
\caption{The cost function $F_{c}$ and the group accuracy rate $\gamma$ of K-means, PSNMF and EP4Orth+ for weak low pass graph signals.}
\label{fig2}
\end{center}
\end{figure*}

We first use the cost function $F_{c}$ and the group accuracy rate $\gamma$ to measure the performance of BE-EEPs, which is \Cref{alg:2}. More specifically, the cost function $F_{c}$ is given by $F_{c}=F(\widehat{\mathbb{C}}, \overline{V})$, where $\widehat{\mathbb{C}}$ represents the EEP we found using either K-means, PSNMF or EP4Orth+, and  $\overline{V}\in\mathbb{R}^{n\times r}$ is a matrix containing the $r$ structural eigenvectors of the Laplacian matrix of the graph under consideration. 
The lower the cost function value, the higher the quality of the partition extracted by the algorithm. Note that each class in the planting partition contains $126$ vertices, denoted as $|C_{i}|=126$, $i=1,\cdots,r$. Let $C^{w}$ represent the sum of the absolute difference between the number of vertices in each class of the output partition and $|C_{i}|$. The group accuracy rate $\gamma$ is given by
$$\gamma=\frac{|V(G)|-C^{w}}{|V(G)|}.$$
Hence, the group accuracy rate of outputs over multiple trials is a value greater than or equal to $0$ and less than or equal to $1$. 
This value indicates the degree of similarity between the number of vertices in the output partitions and the number of vertices in the corresponding classes of the planting partition. The higher the group accuracy, the closer the number of vertices in the output partitions is to the number of vertices in the original partitions.
The two values are shown in \Cref{fig1} and \Cref{fig2}, with the results being the averages over 200 repetitions with different sample sizes.

From \Cref{fig1} and \Cref{fig2}, we observe that for all algorithms, the cost function $F_{c}$ decreases and the group accuracy rate $\gamma$ increases as the number of samples increases. This indicates that, our method BE-EEPs is effective for the blind extraction of EEPs dealing with in both strong and weak low pass graph signals.
Specifically, when extracting from strong low pass graph signals, as shown in the left plot of \Cref{fig1}, a decrease that all three algorithms have in the cost function can be seen for increasing sample size. Obviously, EP4Orth+ yields the highest cost function value relative to PSNMF and K-means, while PSNMF and K-means provide nearly identical cost function values. In the right plot of \Cref{fig1}, we can clearly observe that although PSNMF achieves the highest group accuracy value compared to K-means and EP4Orth+.
When extracting from weak low pass graph signals, as shown in \Cref{fig2}, we observe that there is no significant difference in the cost function $F_{c}$ when comparing the use of K-means, PSNMF and EP4Orth+ for generation. However, it is worth noting that in terms of the group accuracy rate of the generated EEPs, PSNMF slightly outperforms both K-means and EP4Orth+.

\begin{figure}[ht]
	\centering
	\begin{minipage}{1\linewidth}	
		\subfigure[The result for EP4Orth+]{
			\label{fig3-1}
			\includegraphics[width=0.3\linewidth,height=3in]{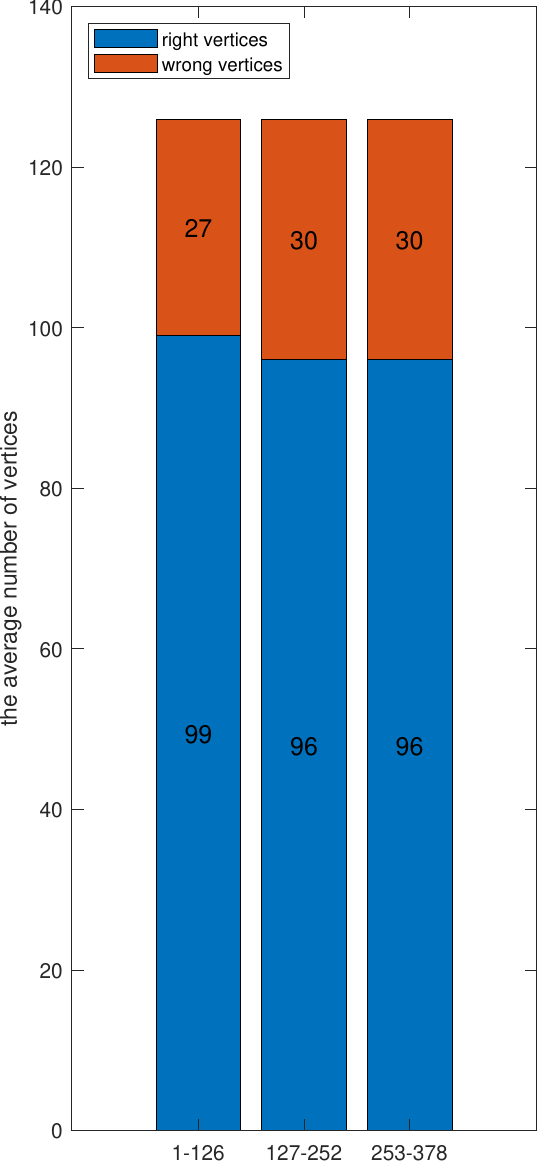}	
		}\noindent
		\subfigure[The result for PSNMF]{
			\label{fig3-2}
			\includegraphics[width=0.3\linewidth,height=3in]{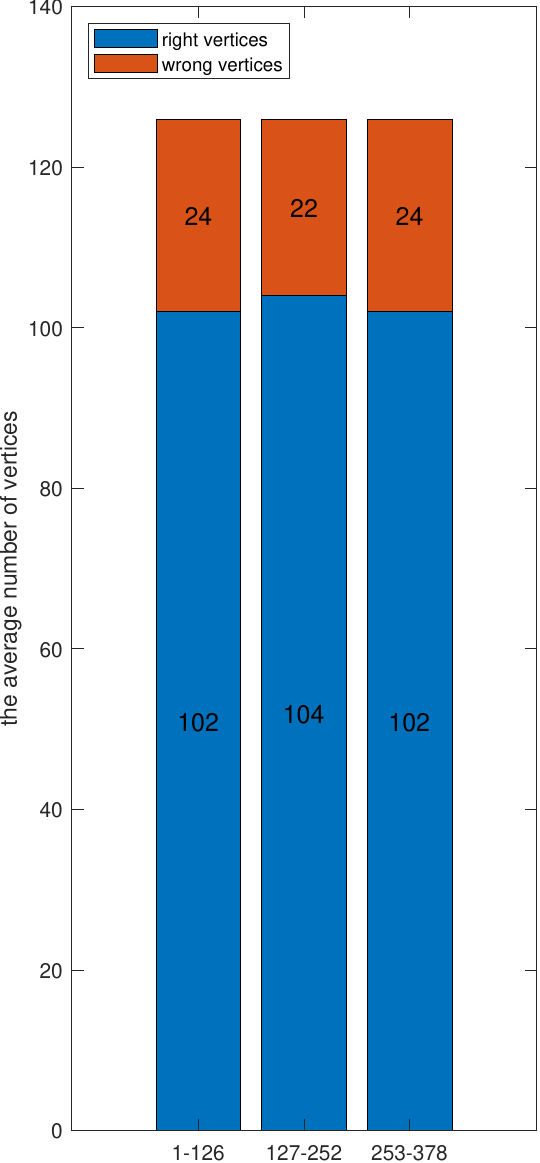}
		}\noindent
		\subfigure[The result for K-means]{
			\label{fig3-3}
			\includegraphics[width=0.3\linewidth,height=3in]{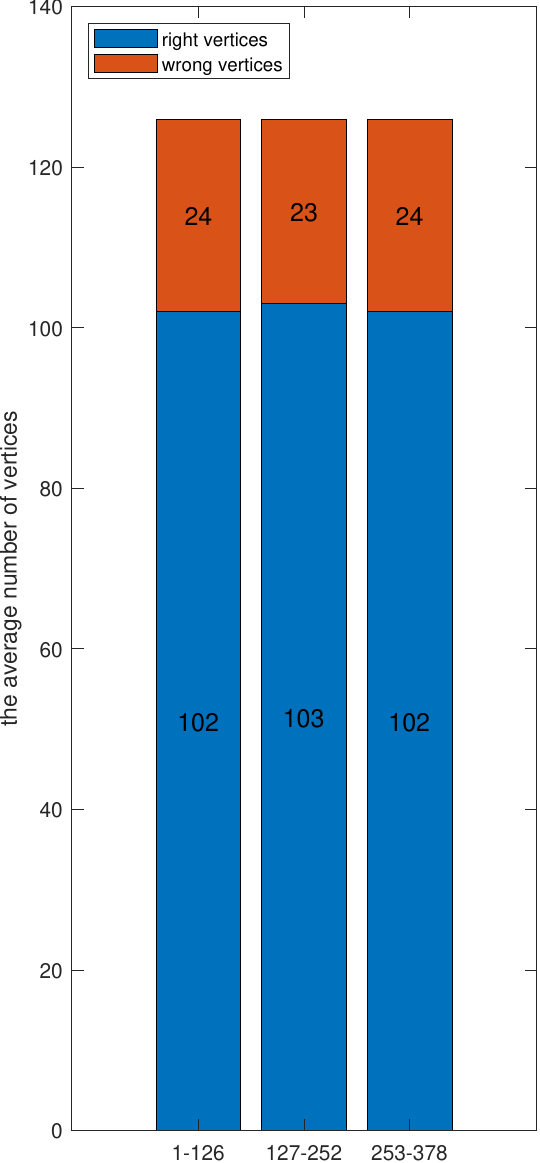}
		}
	\end{minipage}
	\caption{The facet stacked bar chart of EEP extraction from strong low pass graph signals using three methods.}
	\label{fig3}
\end{figure}

\begin{figure}[ht]
	\centering
	\begin{minipage}{1\linewidth}	
		\subfigure[The result for EP4Orth+]{
			\label{fig4-1}
			\includegraphics[width=0.3\linewidth,height=3in]{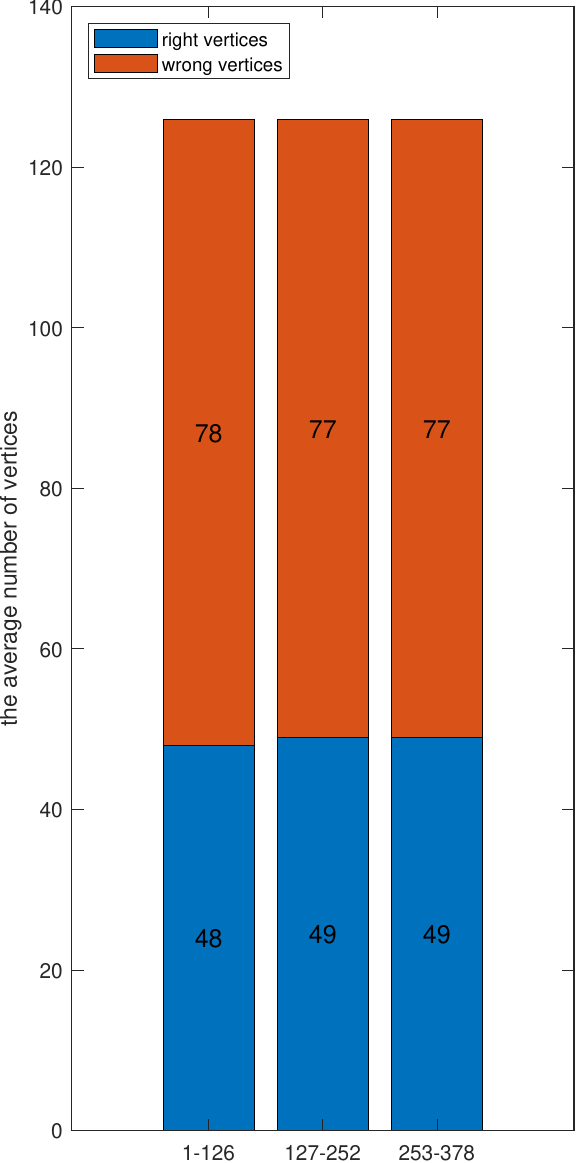}	
		}\noindent
		\subfigure[The result for PSNMF]{
			\label{fig4-2}
			\includegraphics[width=0.3\linewidth,height=3in]{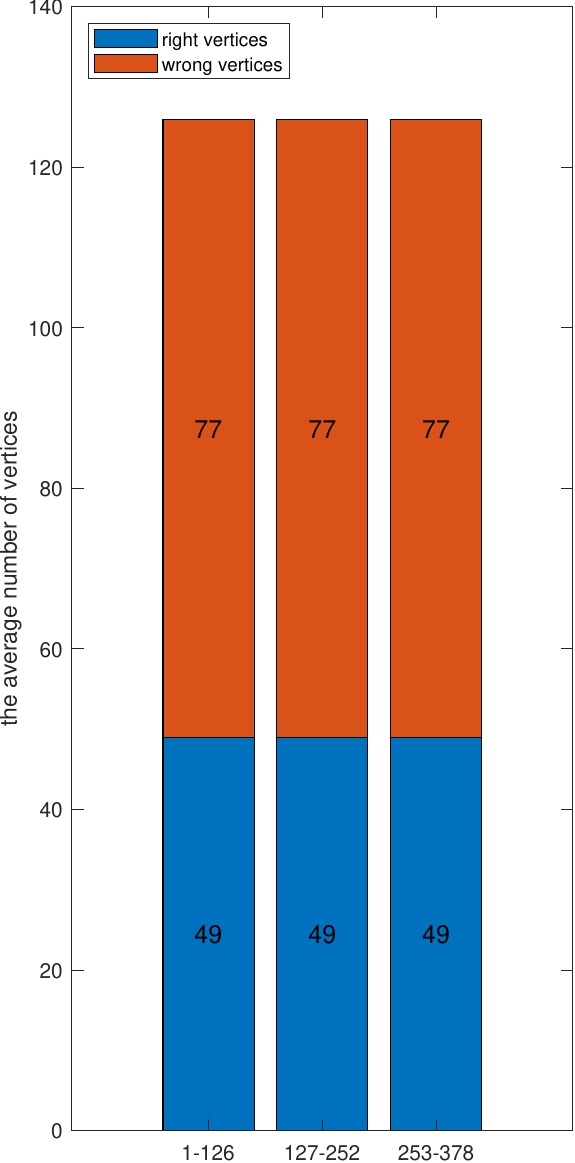}
		}\noindent
		\subfigure[The result for K-means]{
			\label{fig4-3}
			\includegraphics[width=0.3\linewidth,height=3in]{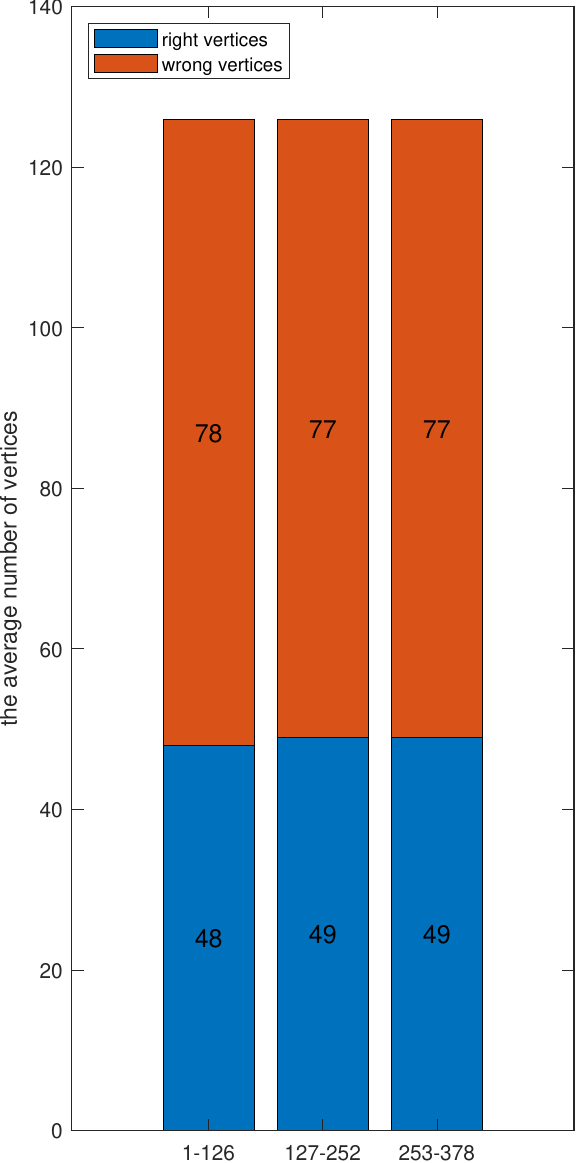}
		}
	\end{minipage}
	\caption{The facet stacked bar chart of EEP extraction from weak low pass graph signals using three methods.}
	\label{fig4}
\end{figure}

To more accurately evaluate the correctness of the partitions extracted by the algorithm, we present \Cref{fig3} and \Cref{fig4}, respectively. Note that we know the true partitions, where vertices $1$-$126$ form one partition, vertices $127$-$252$ form another, and vertices $253$-$378$ form a third partition. In \Cref{fig3} and \Cref{fig4}, each bar represents the number of vertices in a partition. The orange color indicates the number of incorrectly classified vertices, while the blue color represents the number of correctly classified vertices. \Cref{fig3-1}, \Cref{fig3-2}, and \Cref{fig3-3} represent the number of correctly and incorrectly classified vertices under strong low pass graph signals by EP4Orth+, PSNMF, and K-means, respectively. \Cref{fig4-1}, \Cref{fig4-2}, and \Cref{fig4-3} represent the number of correctly and incorrectly classified vertices under weak low pass graph signals by EP4Orth+, PSNMF, and K-means, respectively. The number of vertices in \Cref{fig3} and \Cref{fig4} is the average number of vertices from $181$ experiments, using $300$ samples as an example.

\Cref{fig3} shows that the number of correctly classified vertices provided by PSNMF and K-means is similar and slightly higher than that of EP4Orth+. 
\Cref{fig4} shows that when extracting from weak low pass graph signals, the number of correctly classified vertices generated by the three algorithms is very similar. Therefore, although PSNMF slightly outperforms K-means and EP4Orth+ in terms of group accuracy rate for generating EEPs, when extracting from strong low pass graph signals, K-means and PSNMF perform equally well in both the cost function $F_{c}$ and the number of correctly classified vertices (i.e., low $F_{c}$ and high number of correctly classified vertices), with both slightly better than EP4Orth+. However, when extracting from weak low pass graph signals, $F_{c}$ and the number of correctly partitioned vertices generated by all three algorithms are very similar.

In summary, as mentioned above, our proposed approach BE-EEPs is effective for the blind extraction of EEPs. We find that when extracting EEPs from strong low pass graph signals, PSNMF and K-means perform slightly better than EP4Orth+. However, when extracting EEPs from weak low pass graph signals, all three algorithms are equally effective.

\section{Conclusion}\label{sec7}
In this paper, we addressed the link between low pass graph signals and the extraction of EEPs, alongside addressing the challenge of blind extraction of EEPs under low pass graph filters. To tackle this problem, we proposed a method called BE-EEPs, which reformulated EEP as a matrix optimization model. We introduced the use of K-means, PSNMF and EP4Orth+ to solve this optimization model and offer theoretical support for the proposed method.
The experimental results verified that our method for extracting EEPs from low pass graph signals is effective. In terms of three candidate solvers for solving the resulting matrix optimization model, PSNMF and K-means performed equally well under strong low pass graph signals, while EP4Orth+ was slightly less effective. In more complex weak low pass graph signals, all three algorithms performed equally well.


\vspace{2.5mm}
\noindent {\bf Acknowledgments.} We would like to thank the team of Prof. Jiang Bo from Nanjing Normal University for providing the code for the EP4Orth+ algorithm in \cite{no33}. Without their code, our comparison using EP4Orth+ would not have been possible.

\end{document}